\documentclass[reqno]{amsart}
\usepackage{s-slice}
\usepackage{graphicx}

\let\TSp\thinspace
\def\DSp{\thinspace\thinspace}
\def\lbr{\raise 1pt\hbox{[}}
\def\rbr{\raise 1pt\hbox{]}}

\def\gobble#1{}
\def\tableline{\multispan\numcolumns\hrulefill\cr}
\def\dblhline{\hline height 0.16667em\gobble&\emptyline\cr\hline}

\newbox\tablebox

\newbox\emptybox
\newbox\circlebox
\newbox\graycirclebox
\newbox\squarebox
\newbox\graysquarebox
\newbox\diagbox

\setbox\emptybox\null

\def\PSpreamble#1{%
  gsave /scaleFactor {Resolution 72 div} def
  currentpoint translate #1 scaleFactor mul dup dup dup neg translate scale
  0.4 #1 div setlinewidth
}
\def\PStrailer#1{%
  gsave #1 grestore 0 setgray stroke grestore
}

\def\makecircle#1#2#3{%
  \global\setbox#1\hbox to #2bp{\vbox to #2bp{\vss
    \special{ps: \PSpreamble{#2 2 div} newpath 0 0 1 0 360 arc \PStrailer{#3}}%
  }\hss}%
}
\makecircle\circlebox9{1 setgray eofill}
\makecircle\graycirclebox9{0.75 setgray eofill}

\def\makesquare#1#2#3{%
  \global\setbox#1\hbox to #2bp{\vbox to #2bp{\vss
  }\hss}%
}
\makesquare\squarebox{11}{1 setgray eofill}
\makesquare\graysquarebox{11}{0.75 setgray eofill}

\newdimen\tempPSdima \newdimen\tempPSdimb
\newcount\tempPScnta \newcount\tempPScntb
\def\makediagonal#1#2#3#4{%
  \tempPSdima #2\tempPScnta\tempPSdima
  \tempPSdimb #3\tempPScntb\tempPSdimb
  \ifx#4-\relax
    \def\tempPS{dup 0 exch translate 1 -1 scale}%
  \else
    \def\tempPS{}%
  \fi
  \global\setbox#1\hbox to #2{\vbox to #3{\vss
  }\hss}%
}

\def\zerobox#1{\hbox to 0pt{\hss\vbox to 0pt{\vss #1\vss}\hss}}

\def\joincorners#1#2#3#4#5{%
  \makediagonal\diagbox{#1}{#2}{#5}%
  \vbox{\offinterlineskip
    \copy\diagbox
    \ifx#5-\relax\kern -\ht\diagbox\fi
      \hbox to\wd\diagbox{\zerobox{\copy #3}\hfill}%
    \ifx#5-\relax\kern \ht\diagbox\else\kern -\ht\diagbox\fi
      \hbox to\wd\diagbox{\hfill\zerobox{\copy #4}}%
    \ifx#5-\relax\else\kern \ht\diagbox\fi
  }%
}

\newbox\tablerow
\newbox\KMtable
\newbox\tmpKMbox

\newcount\rowcount

\newdimen\rowheight
\newdimen\botheight
\newdimen\lastentrywidth
\newdimen\leftKMskip
\newdimen\rightKMskip
\newdimen\KMwidth
\newdimen\entrymargin	
\newdimen\Tentrymargin	

\def\showKMpairs#1{%
  \edef\KMpairs{#1}%
  \setbox\tmpKMbox\vbox{}%
  \setbox\tablebox\vbox{%
    \offinterlineskip
    \unvbox\tablebox
    \rowcount 4\setbox\tablebox\vbox{}%
    \loop
      \advance\rowcount -1\relax
      \global\setbox\tablerow\lastbox\unskip
      \global\setbox\tmpKMbox\vbox{\copy\tablerow\unvbox\tmpKMbox}%
    \ifnum\rowcount>0\repeat
    \unvcopy\tmpKMbox
  }%

  \rowheight.4pt\advance\rowheight\ht\tablerow \advance\rowheight\dp\tablerow
  \botheight-\rowheight\advance\botheight.2pt
  \advance\botheight\ht\tmpKMbox \advance\botheight\dp\tmpKMbox

  \setbox\tmpKMbox\hbox{\DSp}\entrymargin\wd\tmpKMbox\Tentrymargin\entrymargin
  \setbox\tmpKMbox\hbox{1\kern.4pt}\advance\entrymargin .5\wd\tmpKMbox
  \setbox\tmpKMbox\hbox{\lbr1\rbr\kern.2pt}\advance\Tentrymargin .5\wd\tmpKMbox

  \setbox\KMtable\hbox{\vrule width0pt height\ht\tablebox depth\dp\tablebox}%
  \hbox{%
    \global\lastentrywidth 0pt\relax
    \unhbox\tablerow\unskip
    \loop
      \setbox\tmpKMbox\lastbox\unskip
    \ifhbox\tmpKMbox
      \global\setbox\KMtable\hbox{%
        \hbox to\wd\tmpKMbox{\expandafter\processKHpair\KMpairs\theend}%
        \unhbox\KMtable
      }%
      \global\lastentrywidth\wd\tmpKMbox
    \repeat
  }%
}

\def\empty{}
\def\processKHpair#1,#2:#3\theend{%
  \gdef\KMpairs{#3}%
  \def\firstargument{#1}%
  \ifx\firstargument\empty\else
    \vbox{\offinterlineskip
      \hbox{\csname KM#1\endcsname}%
      \kern #2.5\rowheight\kern\botheight\relax
    }\hss
  \fi
}

\def\setKMsizes{%
  \KMwidth\wd\tmpKMbox \advance\KMwidth\lastentrywidth
  \advance\KMwidth -\leftKMskip \advance\KMwidth -\rightKMskip

  \kern\leftKMskip\kern -.2pt
}

\def\simplepair{%
  \setKMsizes
  \joincorners\KMwidth{2\rowheight}\circlebox\circlebox{}%
}

\def\doublepair{%
  \vbox{\offinterlineskip
    \hbox{\simplepair}\kern -\rowheight\hbox{\simplepair}%
  }%
}

\def\specialVpair{%
  \advance\lastentrywidth\wd\tmpKMbox\setKMsizes
  \joincorners\KMwidth\rowheight\graycirclebox\graycirclebox{}%
}

\def\torsionpair{%
  \vbox{\offinterlineskip
    \hbox{%
      \leftKMskip\wd\tmpKMbox \advance\leftKMskip .5\lastentrywidth
      \setKMsizes\joincorners\KMwidth\rowheight\circlebox\emptybox-%
    }%
    \hbox{\setKMsizes\joincorners\KMwidth\rowheight\circlebox\graysquarebox{}}%
  }%
}

\long\def\tbl#1#2{$$\hbox{#2}$$\caption{#1}}

\begin{document}
\title[Rasmussen invariant and sliceness of knots]{Rasmussen invariant,\\
slice-Bennequin inequality,\\ and sliceness of knots}
\author[A.~Shumakovitch]{Alexander N.~Shumakovitch}
\address{Department of Mathematics, The George Washington University,
Phillips Hall, 801\ 22nd St. NW, Suite \#739, Washington, DC 20052, U.S.A.}
\email{Shurik@gwu.edu}

\begin{abstract}
We use recently introduced Rasmussen invariant to find knots that are
topologically locally-flatly slice but not smoothly slice. We note that this
invariant can be used to give a combinatorial proof of the slice-Bennequin
inequality. Finally, we compute the Rasmussen invariant for quasipositive
knots and show that most of our examples of non-slice knots are not
quasipositive and, to the best of our knowledge, were previously unknown.
\end{abstract}

\maketitle

\section{Rasmussen invariant and the slice-Bennequin inequality}
In~\cite{Rasmussen-s} Jacob Rasmussen used the theory of knot (co)homology
developed by Mikhail Khovanov~\cite{Khovanov-categor} and results of Eun
Soo Lee~\cite{Lee-alt-links} to introduce a new invariant $s$ of knots in
$S^3$. This invariant takes values in even integers. Its main properties are
summarized as follows.

\begin{theorem}[{Rasmussen~\cite[Theorems~1--4]{Rasmussen-s}}]
\label{thm:Rasmussen-main}
Let $K$ be a knot in $S^3$. Then
\begin{enumerate}
\item $s$ gives a lower bound on the slice ($4$-dimensional) genus $g_s(K)$
of $K$:
\begin{equation}\label{eq:s-g-bound}
|s(K)|\le2g_s(K);
\end{equation}
\item $s$ induces a homomorphism from $Conc(S^3)$, the concordance
group of knots in $S^3$, to $\Z$;
\item If $K$ is alternating, then $s(K)=\Gs(K)$, where $\Gs(K)$ is the
classical knot signature of $K$;
\item If $K$ can be represented by a positive diagram $D$, then
\begin{equation}\label{eq:s-positive-knots}
s(K)=2g_s(K)=2g(K)=n(D)-O(D)+1,
\end{equation}
where $n(D)$ and $O(D)$ are the number of crossings and Seifert circles of
$D$, respectively, and $g(K)$ is the ordinary (3-dimensional) genus of $K$.
\end{enumerate}
\end{theorem}

\begin{corollary}[{\cite[Corollary~4.3]{Rasmussen-s}}]\label{cor:s-Xing-change}
Let $K_+$ and $K_-$ be two knots that are different at a single crossing that
is positive in $K_+$ and negative in $K_-$. Then
\begin{equation}\label{eq:s-change-Xing}
s(K_-)\le s(K_+)\le s(K_-)+2.
\end{equation}
\end{corollary}

Equality~\eqref{eq:s-positive-knots} can be easily generalized to arbitrary
knots. It becomes an inequality then.

\begin{lemma}\label{lem:s-inequality}
Let $K$ be a knot represented by an (oriented) diagram $D$ and let $w(D)$ be
the {\em writhe number} of $D$, that is, the number of positive crossings of
$D$ minus the number of negative ones:
$\displaystyle w(D)=\#\left\{
\vcenter{\hbox{\includegraphics[scale=0.55]{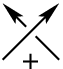}}}\right\}-\#\left\{
\vcenter{\hbox{\includegraphics[scale=0.55]{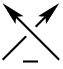}}}\right\}$.
Then
\begin{equation}\label{eq:s-writhe-Seifert}
s(K)\ge w(D)-O(D)+1.
\end{equation}
\end{lemma}
\begin{proof}
Use induction on the number of negative crossings in $D$. If all crossings of
$D$ are positive, the Lemma follows from~\eqref{eq:s-positive-knots}. The
inequality~\eqref{eq:s-writhe-Seifert} is an equality then. When a positive
crossing of $D$ is changed into a negative one, the right hand side
of~\eqref{eq:s-writhe-Seifert} decreases by $2$, while the left hand side
decreases by at most $2$ because of~\eqref{eq:s-change-Xing}. Hence, the
inequality is preserved.
\end{proof}

Lemma~\ref{lem:s-inequality} implies the slice-Bennequin inequality that
was originally proved by Rudolph~\cite{Rudolph-quasipositivity} using the
gauge theory. The theory developed by Khovanov, Lee, and Rasmussen provides
the first purely combinatorial proof of this inequality.

\begin{corollary}[Slice-Bennequin Inequality, cf. \cite{Rudolph-quasipositivity}]
\label{cor:slice-bennequin}
Let $\Gb$ be a braid with $k$ strands and let $\widehat\Gb$ be its closure.
Denote by $\chi_s(\widehat\Gb)$ the greatest Euler characteristic of an
oriented surface (without closed components) smoothly embedded in $D^4$ with
boundary $\widehat\Gb$. Then
\begin{equation}\label{eq:slice-bennequin}
\chi_s(\widehat\Gb)\le k-w(\Gb).
\end{equation}
\end{corollary}
\begin{proof}
If $\widehat\Gb$ is a knot, then~\eqref{eq:s-writhe-Seifert}
and~\eqref{eq:s-g-bound} imply that
\begin{equation}
g_s(\widehat\Gb)\ge \frac{w(\Gb)-k+1}2.
\end{equation}
It remains to notice that $\chi_s=1-2g_s$ for knots.

Assume now that $\widehat\Gb$ is a link. Let $\Gb^+$ be a braid obtained from
$\Gb$ by removal from the braid word representing $\Gb$ of all the standard
generators that appear with negative exponents. For example, if
$\Gb=\Gs_2\Gs_1^{-1}\Gs_2$, then $\Gb^+=\Gs_2^2$. Inserting a cancelling pair
of generators $\Gs_i\Gs_i^{-1}$ into $\Gb$ changes neither $w(\Gb)$ nor
$\chi_s(\widehat\Gb)$, but adds a crossing to $\Gb^+$, so one can assume
without a loss of generality that the closure $\widehat{\Gb^+}$ of $\Gb^+$ is
a knot.

Since $\widehat{\Gb^+}$ is a (positive) knot,~\eqref{eq:slice-bennequin} holds
true for it (in fact, it is an equality). Now, addition of a negative crossing
to a braid increases the right-hand side of~\eqref{eq:slice-bennequin} by
exactly $1$. On the other hand, the following Lemma shows that $\chi_s$ can
not change by more than $1$. This completes the proof.
\end{proof}

\begin{lemma}
Let $\Gb$ and $\Gb'$ be two braids such that $\Gb=w_1w_2$ and
$\Gb'=w_1\Gs_i^\varepsilon w_2$, where $w_1$ and $w_2$ are some braid words,
$\Gs_i$ is a standard braid group generator, and $\varepsilon=\pm1$. Let
$\widehat\Gb$ and $\widehat{\Gb'}$ be the corresponding closures. Then
$|\chi_s(\widehat\Gb)-\chi_s(\widehat{\Gb'})|\le1$.
\end{lemma}
\begin{proof}
Let $S$ be an oriented surface (without closed components) smoothly embedded
in $D^4$ with $\partial S=\widehat\Gb$ and $\chi(S)=\chi_s(\widehat\Gb)$.
Addition of a twisted band to $S$ at the place where a crossing is added to
$\Gb$ produces a smoothly embedded surface $S'$ with $\partial
S'=\widehat{\Gb'}$. Then $\chi(S')=\chi(S)-1$ and
$\chi_s(\widehat{\Gb'})\ge\chi_s(\widehat\Gb)-1$. On the other hand,
$\Gb=w_1\Gs_i^{-\varepsilon}\Gs_i^\varepsilon w_2$. Then
$\Gb=w_1\Gs_i^{-\varepsilon}w_2'$ and $\Gb'=w_1w_2'$ with
$w_2'=\Gs_i^\varepsilon w_2$. Repeating the previous argument, one obtains
that $\chi_s(\widehat\Gb)\ge\chi_s(\widehat{\Gb'})-1$.
\end{proof}

The slice-Bennequin inequality leads to a formula for the Rasmussen invariant
of (strongly) quasipositive knots. We use this formula in
section~\ref{sec:sliceness}.

\begin{definitions}
1.~A knot $K$ is said to be {\em quasipositive} if it is the closure of a braid
that has the form $(w_1\Gs_{j_1}w_1^{-1})(w_2\Gs_{j_2}w_2^{-1})\cdots
(w_p\Gs_{j_p}w_p^{-1})$, where $\Gs_i$ are the standard generators of the
corresponding braid group, and $w_i$ are braid words.

\noindent
2.~A knot $K$ is said to be {\em strongly quasipositive} if it is the closure
of a braid that has the form
$\Gs_{{i_1,j_1}}\Gs_{{i_2,j_2}}\cdots\Gs_{{i_p,j_p}}$, where
$\Gs_{i,j}=(\Gs_i\Gs_{i+1}\cdots\Gs_{j-2})\Gs_{j-1}
(\Gs_i\Gs_{i+1}\cdots\Gs_{j-2})^{-1}$ for $j\ge i+2$ and
$\Gs_{i,i+1}=\Gs_i$.
\end{definitions}

\begin{figure}
\centerline{\hbox{%
\vbox{\halign{\hfil#\hfil\cr
\includegraphics[scale=0.50]{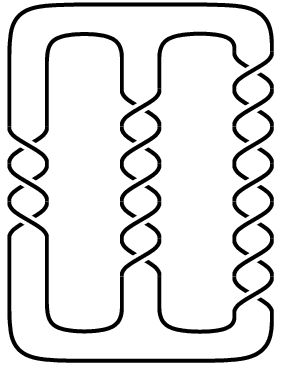}\cr a.\cr}}%
\qquad\qquad
\vbox{\halign{\hfil#\hfil\cr
\includegraphics[scale=0.60]{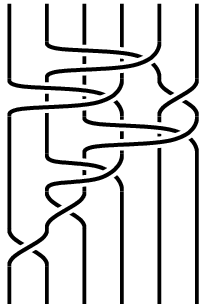}\cr b.\cr}}$
\qquad\qquad
\vbox{\halign{\hfil#\hfil\cr
\includegraphics[scale=0.60]{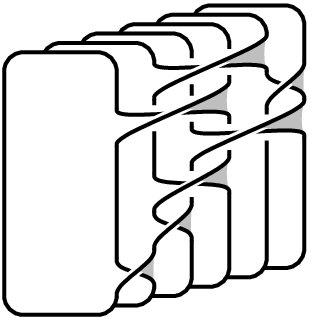}\cr c.\cr}}$
}}
\caption{The $(-3,5,7)$-pretzel knot, its representation as the closure of a
strongly quasipositive braid
$\Gs_1\,\Gs_2\,\Gs_{2,4}\,\Gs_{3,6}\,\Gs_{1,4}\,\Gs_5\,\Gs_{2,5}$, and
the corresponding Seifert surface.} \label{fig:3-5-7-pretzel}
\end{figure}

The $(-3,5,7)$-pretzel knot depicted in Figure~\ref{fig:3-5-7-pretzel}.a is a
strongly quasipositive knot. It is the closure of the braid
$\Gs_1\,\Gs_2\,\Gs_{2,4}\,\Gs_{3,6}\,\Gs_{1,4}\,\Gs_5\,\Gs_{2,5}$ (see
Figure~\ref{fig:3-5-7-pretzel}.b).

\begin{proposition}\label{prop:s-for-quasipositive}
Let $K$ be a knot that can be represented as the closure of a quasipositive
braid $\Gb$ with $k$ strands and $b$ bands, that is, $\Gb$ is a product of $b$
factors of the form $w\Gs_i w^{-1}$. Then
\begin{equation}\label{eq:s-for-quasipositive}
s(K)=2g_s(K)=b-k+1.
\end{equation}
If, moreover, $\Gb$ is strongly quasipositive, then
\begin{equation}\label{eq:s-for-strongly-quasipositive}
s(K)=2g(K)=2g_s(K)=b-k+1.
\end{equation}
\end{proposition}
\begin{proof}
It is clear that $w(K)=b$ and $O(K)=k$. Hence,
\begin{equation*}
b-k+1\le s(K)\le|s(K)|\le2g_s(K)\le2g(K)
\end{equation*}
by~\eqref{eq:s-writhe-Seifert} and~\eqref{eq:s-g-bound}. On the other hand,
for a (strongly) quasipositive knot one can explicitly construct a surface $S$
smoothly embedded in $D^4$ (respectively, $S^3$), such that $\partial S=K$ and
the Euler characteristic $\chi(S)$ of $S$ is $k-b$ (see
Figure~\ref{fig:3-5-7-pretzel}.c that illustrates this construction in the
strongly quasipositive case). Since $S$ has a single boundary component, its
genus equals $\frac{1-\chi(S)}2$. It follows that $g_s(K)$ (respectively,
$g(K)$) does not exceed $\frac{b-k+1}2$. This finishes the proof.
\end{proof}

\begin{remark}
One of the main applications of the $s$-invariant provided
in~\cite{Rasmussen-s} was a purely combinatorial proof of the Milnor
Conjecture~\cite{Milnor-conjecture} that was originally proved by Kronheimer
and Mrowka~\cite{Kronheimer-Mrowka} using the gauge theory. More specifically,
Rasmussen showed that the slice genus of a $(p,q)$-torus knot with $p,q>0$
is $(p-1)(q-1)/2$. In fact, the original question posed by Milnor and answered
in~\cite{Kronheimer-Mrowka} is more general
(see~\cite[Remark~10.9]{Milnor-conjecture}). It asks whether the unknotting
number $u(K)$ of an algebraic knot $K$ equals its genus $g(K)$. Here, an {\sl
algebraic knot} is the knot associated to an isolated singular point of a
complex algebraic curve in $\C^2$ by intersecting it with a $3$-dimensional
sphere of a sufficiently small radius centered at the singularity. For
example, a $(p,q)$-torus knot corresponds to the singular curve $z^p+w^q=0$.

It is well-known that $g_s(K)\le u(K)\le g(K)$ for any algebraic knot $K$ and
that it can be represented by a positive diagram obtained as a closure of a
positive braid. The general Milnor Conjecture now follows from
\eqref{eq:s-positive-knots} straightforwardly. I am grateful to Sergei Chmutov
for pointing out this fact to me.
\end{remark}

\begin{remark}
Charles Livingston~\cite{Livingston-tau} used the Ozsv\'ath-Szab\'o knot
invariant $\Gt$~\cite{Ozsvath-Szabo-tau} to give new proofs to several
results of Lee Rudolph~\cite{Rudolph-quasipositivity, Rudolph-sliceness} on
the slice genus, including the slice-Bennequin Inequality. Invariants $s$ and
$\Gt$ share many of their main properties and our approach is similar to the
Livingston's one. The key difference is that the Rasmussen invariant is
defined combinatorially, while the Ozsv\'ath-Szab\'o invariant is based on the
theory of knot Floer homology. It was originally conjectured that
$s(K)=2\Gt(K)$ for every knot $K$, but counter-examples were later
found~\cite{Hedden-Ording,Livingston-slice}.
\end{remark}

\begin{remark}
Relation between the Rasmussen invariant and the slice-Bennequin
inequality was independently observed by several other authors, including
Olga Pla\-me\-nevskaya~\cite{Plamenevskaya-transverse} and Alexander
Stoimenow~\cite{Stoimenow-sliceness}. After the original version of this
paper was published, Tomomi Kawamura~\cite{Kawamura} used the Rasmussen
invariant to prove a sharper slice-Bennequin inequality.
\end{remark}

\section{Sliceness of knots}\label{sec:experiments}
\label{sec:sliceness}

In many cases one can easily compute $s(K)$ from the Khovanov homology of $K$.
For a given knot $K$, let $h^{i,j}(K)=\dim_\Q(\CalH^{i,j}(K)\otimes\Q)$ be the
ranks of its homology and $Kh(K)(t,q)=\sum_{i,j}t^iq^jh^{i,j}(K)$ be the
corresponding Poincar\'e polynomial in variables $t$ and $q$. Denote by
$hw(K)$ the {\em homological width} of $K$, that is, the minimal number of
adjacent diagonals $2i-j=const$ that support the homology of $K$.

In was shown by Rasmussen~\cite[Proposition 5.2]{Rasmussen-s} that for all
knots $K$ with $hw(K)\le 3$, one has
\begin{equation}\label{eq:s-equation}
Kh(K)=q^{s(K)-1}\big(1+q^2+(1+tq^4)Kh'(K)\big),
\end{equation}
where $Kh'(K)$ is some (Laurent) polynomial in $t$ and $q$ with non-negative
coefficients. In fact, Rasmussen's arguments can be applied to a more general
case.
\begin{statement}\label{prop:s-from-Kh}
Let $K$ be a knot. Assume that $h^{i,j}(K)h^{i+1,j+4(n-1)}(K)=0$ for all $i$,
$j$, and $n\ge3$ (this is automatically the case if $hw(K)\le 3$).
Then~\eqref{eq:s-equation} holds true for $K$.
\end{statement}
\begin{proof}
Construction of the Rasmussen invariant is based on a spectral sequence
structure on the Khovanov chain complex that is due to
Lee~\cite{Lee-alt-links}. The differential $d_n$ in this spectral sequence
has bidegree $(1, 4(n-1))$. The condition on $h^{i,j}$ implies that $d_n$ is
trivial for all $n\ge 3$. The rest of the arguments is the same as
in~\cite[Proposition 5.2]{Rasmussen-s}.
\end{proof}

It is possible for a knot to have homological width $4$, but still satisfy the
condition of~\ref{prop:s-from-Kh} (see Table~\ref{tbl:16n_809057-knot}). On
the other hand, the knot $16^n_{864894}$\footnote{We enumerate knots according
to the convention from Knotscape~\cite{Knotscape}, due to Hoste and
Thistlethwaite. For example, the knot $13^n_{1496}$ is a non-alternating knot
number 1496 with 13 crossings.} may theoretically have $d_3^{-1,-7}\not=0$,
since $h^{-1,-7}=h^{0,1}=1$ (see Table~\ref{tbl:16n_864894-knot}). Hence, its
Rasmussen invariant can equal either $0$ or $-2$. Let us demonstrate that it
is indeed the former.

\begin{table}[t]
\tbl{Homology ranks of the knot $16^n_{809057}$. The homological width is
4, but the differential $d_3$ of bidegree $(1,8)$ is trivial.}
{\def\DSp{\kern 0.2em}%
\def\emptyline{&&&&&&&&&&&&&&&&&}%
\def\numcolumns{19}%

\setbox\tablebox\vbox{\small\let\hline\tableline\offinterlineskip\ialign{%
\vrule\TSp\vrule #\relax&\DSp\hfil #\DSp\vrule\TSp\vrule&
\DSp\hfil\bf #\hfil\DSp\vrule&
\DSp\hfil\bf #\hfil\DSp\vrule&
\DSp\hfil\bf #\hfil\DSp\vrule&
\DSp\hfil\bf #\hfil\DSp\vrule&
\DSp\hfil\bf #\hfil\DSp\vrule&
\DSp\hfil\bf #\hfil\DSp\vrule&
\DSp\hfil\bf #\hfil\DSp\vrule&
\DSp\hfil\bf #\hfil\DSp\vrule&
\DSp\hfil\bf #\hfil\DSp\vrule&
\DSp\hfil\bf #\hfil\DSp\vrule&
\DSp\hfil\bf #\hfil\DSp\vrule&
\DSp\hfil\bf #\hfil\DSp\vrule&
\DSp\hfil\bf #\hfil\DSp\vrule&
\DSp\hfil\bf #\hfil\DSp\vrule&
\DSp\hfil\bf #\hfil\DSp\vrule&
\DSp\hfil\bf #\hfil\DSp\vrule&#\TSp\vrule\cr
\dblhline
height 11pt depth 4pt&&
\rm\DSp-9\DSp&
\rm\DSp-8\DSp&
\rm\DSp-7\DSp&
\rm\DSp-6\DSp&
\rm\DSp-5\DSp&
\rm\DSp-4\DSp&
\rm\DSp-3\DSp&
\rm\DSp-2\DSp&
\rm\DSp-1\DSp&
\rm\DSp0\DSp&
\rm\DSp1\DSp&
\rm\DSp2\DSp&
\rm\DSp3\DSp&
\rm\DSp4\DSp&
\rm\DSp5\DSp&
\rm\DSp6\DSp&\cr\dblhline
height 10pt depth 3pt&9&
&
&
&
&
&
&
&
&
&
&
&
&
&
&
&
\DSp1\DSp&
\cr\hline
height 10pt depth 3pt&7&
&
&
&
&
&
&
&
&
&
&
&
&
&
&
\DSp2\DSp&
&
\cr\hline
height 10pt depth 3pt&5&
&
&
&
&
&
&
&
&
&
&
&
&
&
\DSp2\DSp&
1&
&
\cr\hline
height 10pt depth 3pt&3&
&
&
&
&
&
&
&
&
&
&
&
\DSp2\DSp&
\DSp3\DSp&
2&
&
&
\cr\hline
height 10pt depth 3pt&1&
&
&
&
&
&
&
&
&
&
&
\DSp2\DSp&
3&
2&
&
&
&
\cr\hline
height 10pt depth 3pt&-1&
&
&
&
&
&
&
&
&
&
\DSp4\DSp&
4&
3&
&
&
&
&
\cr\hline
height 10pt depth 3pt&-3&
&
&
&
&
&
&
&
&
4&
5&
3&
&
&
&
&
&
\cr\hline
height 10pt depth 3pt&-5&
&
&
&
&
&
&
&
3&
4&
2&
&
&
&
&
&
&
\cr\hline
height 10pt depth 3pt&-7&
&
&
&
&
&
\DSp1\DSp&
3&
4&
2&
&
&
&
&
&
&
&
\cr\hline
height 10pt depth 3pt&-9&
&
&
&
&
\DSp1\DSp&
2&
3&
\DSp1\DSp&
&
&
&
&
&
&
&
&
\cr\hline
height 10pt depth 3pt&-11&
&
&
&
&
1&
3&
&
&
&
&
&
&
&
&
&
&
\cr\hline
height 10pt depth 3pt&-13&
&
&
\DSp1\DSp&
1&
2&
&
&
&
&
&
&
&
&
&
&
&
\cr\hline
height 10pt depth 3pt&-15&
&
\DSp1\DSp&
&
&
&
&
&
&
&
&
&
&
&
&
&
&
\cr\hline
height 10pt depth 3pt&-17&
&
1&
&
&
&
&
&
&
&
&
&
&
&
&
&
&
\cr\hline
height 10pt depth 3pt&-19&
\DSp1\DSp&
&
&
&
&
&
&
&
&
&
&
&
&
&
&
&
\cr
\dblhline
}}

\box\tablebox
}
\label{tbl:16n_809057-knot}
\end{table}

\begin{table}[t]
\tbl{Homology ranks of the knot $16^n_{864894}$. Encircled entries show
that $d_3^{-1,-7}$ can possibly be non-trivial.}
{\def\DSp{\kern 0.2em}%
\def\emptyline{&&&&&&&&&&&&&&&&&}%
\def\numcolumns{19}%

\setbox\tablebox\vbox{\small\let\hline\tableline\offinterlineskip\ialign{%
\vrule\TSp\vrule #\relax&\DSp\hfil #\DSp\vrule\TSp\vrule&
\DSp\hfil\bf #\hfil\DSp\vrule&
\DSp\hfil\bf #\hfil\DSp\vrule&
\DSp\hfil\bf #\hfil\DSp\vrule&
\DSp\hfil\bf #\hfil\DSp\vrule&
\DSp\hfil\bf #\hfil\DSp\vrule&
\DSp\hfil\bf #\hfil\DSp\vrule&
\DSp\hfil\bf #\hfil\DSp\vrule&
\DSp\hfil\bf #\hfil\DSp\vrule&
\DSp\hfil\bf #\hfil\DSp\vrule&
\DSp\hfil\bf #\hfil\DSp\vrule&
\DSp\hfil\bf #\hfil\DSp\vrule&
\DSp\hfil\bf #\hfil\DSp\vrule&
\DSp\hfil\bf #\hfil\DSp\vrule&
\DSp\hfil\bf #\hfil\DSp\vrule&
\DSp\hfil\bf #\hfil\DSp\vrule&
\DSp\hfil\bf #\hfil\DSp\vrule&#\TSp\vrule\cr
\dblhline
height 11pt depth 4pt&&
\rm\DSp-9\DSp&
\rm\DSp-8\DSp&
\rm\DSp-7\DSp&
\rm\DSp-6\DSp&
\rm\DSp-5\DSp&
\rm\DSp-4\DSp&
\rm\DSp-3\DSp&
\rm\DSp-2\DSp&
\rm\DSp-1\DSp&
\rm\DSp0\DSp&
\rm\DSp1\DSp&
\rm\DSp2\DSp&
\rm\DSp3\DSp&
\rm\DSp4\DSp&
\rm\DSp5\DSp&
\rm\DSp6\DSp&\cr\dblhline
height 10pt depth 3pt&9&
&
&
&
&
&
&
&
&
&
&
&
&
&
&
&
\DSp1\DSp&
\cr\hline
height 10pt depth 3pt&7&
&
&
&
&
&
&
&
&
&
&
&
&
&
&
\DSp1\DSp&
&
\cr\hline
height 10pt depth 3pt&5&
&
&
&
&
&
&
&
&
&
&
&
&
&
\DSp1\DSp&
1&
&
\cr\hline
height 10pt depth 3pt&3&
&
&
&
&
&
&
&
&
&
&
&
\DSp2\DSp&
\DSp2\DSp&
1&
&
&
\cr\hline
height 10pt depth 3pt&1&
&
&
&
&
&
&
&
&
&
\!{\rm\LARGE\textcircled{\normalsize\bf 1}}\!&
\DSp2\DSp&
1&
1&
&
&
&
\cr\hline
height 10pt depth 3pt&-1&
&
&
&
&
&
&
&
&
&
\DSp3\DSp&
3&
2&
&
&
&
&
\cr\hline
height 10pt depth 3pt&-3&
&
&
&
&
&
&
&
\DSp1\DSp&
\DSp4\DSp&
3&
1&
&
&
&
&
&
\cr\hline
height 10pt depth 3pt&-5&
&
&
&
&
&
&
\DSp1\DSp&
2&
2&
1&
&
&
&
&
&
&
\cr\hline
height 10pt depth 3pt&-7&
&
&
&
&
&
\DSp1\DSp&
3&
4&
\!\!{\rm\LARGE\textcircled{\normalsize\bf 1}}\!\!&
&
&
&
&
&
&
&
\cr\hline
height 10pt depth 3pt&-9&
&
&
&
&
\DSp2\DSp&
3&
2&
&
&
&
&
&
&
&
&
&
\cr\hline
height 10pt depth 3pt&-11&
&
&
&
\DSp1\DSp&
1&
2&
&
&
&
&
&
&
&
&
&
&
\cr\hline
height 10pt depth 3pt&-13&
&
&
\DSp1\DSp&
2&
\DSp2\DSp&
&
&
&
&
&
&
&
&
&
&
&
\cr\hline
height 10pt depth 3pt&-15&
&
\DSp1\DSp&
1&
&
&
&
&
&
&
&
&
&
&
&
&
&
\cr\hline
height 10pt depth 3pt&-17&
&
1&
&
&
&
&
&
&
&
&
&
&
&
&
&
&
\cr\hline
height 10pt depth 3pt&-19&
\DSp1\DSp&
&
&
&
&
&
&
&
&
&
&
&
&
&
&
&
\cr
\dblhline
}}

\box\tablebox
}
\label{tbl:16n_864894-knot}
\end{table}

For a given knot $K$, denote by $rk(K)$ and $\widetilde{rk}(K)$ the total
ranks of its standard and {\em reduced} (see~\cite{Khovanov-patterns} for the
definition) Khovanov homologies, respectively. In particular,
$rk(K)=\sum_{i,j}h^{i,j}(K)=Kh(K)(1,1)$.

\begin{lemma}[{\cite{Khovanov-Frobenius}, page 189}]\label{lem:rank-difference}
Let $K$ be a knot. If $rk(K)-\widetilde{rk}(K)=1$, then~\eqref{eq:s-equation}
holds true for $K$.
\end{lemma}

As it turns out, $rk(16^n_{864894})=66$ and $\widetilde{rk}(16^n_{864894})=65$.
Lemma~\ref{lem:rank-difference} implies that \eqref{eq:s-equation} holds true
for this knot and its Rasmussen invariant can be computed to be equal $0$.

\begin{remark}
Bar-Natan's conjecture~\cite{BarNatan-categor} that~\eqref{eq:s-equation}
holds true for every knot for some integer $s\{K\}$ is still open. Moreover,
no examples are known with $s\{K\}$ different from $s(K)$.
\end{remark}

Given a knot $K$, denote by $\GD_K(t)$ its Alexander polynomial. M.~Freedman
proved~\cite{Freedman-concordance} that if $\GD_K(t)=1$, then $K$ is
topologically locally-flatly slice. We used Knotscape~\cite{Knotscape} to list
all the non-alternating knots with up to 16 crossings that have $\GD=1$. In
total, there are $699$ such knots (not counting the mirror images). The first
one of them has 11 crossings (see Table~\ref{tbl:num-non-slice}). We did not
consider alternating knots, since the Rasmussen invariant equals the signature
for them (see Theorem~\ref{thm:Rasmussen-main}) and the signature is known to
be $0$ when $\GD=1$.

\begin{table}[ht]
\tbl{Number of non-alternating knots that have Alexander polynomial 1 and
those among them that have non-zero Rasmussen invariant.}
{\normalsize\def\arraystretch{1.2}%
\begin{tabular}{|l||*{6}{c|}}
\hline
Number of crossings&11&12&13&14&15&16 \\
\hline\hline
Number of non-alternating knots&185&888&5110&27436&168030&1008906 \\
\hline
\quad among them with $\GD=1$&2&2&15&36&145&499 \\
\hline
\qquad among them with $s\not=0$&0&0&1&1&15&65 \\
\hline
\end{tabular}
}
\label{tbl:num-non-slice}
\end{table}

Next, we used {\tt KhoHo}, a program for computing and studying Khovanov
homology~\cite{Me-KhoHo}, to find the homology of all the knots with
$\GD=1$. Fortunately for us, most of the knots considered (in particular, all
knots with at most $15$ crossings) have homological width $3$. There are 42
knots with 16 crossings and homological width $4$ that satisfy the condition
of~\ref{prop:s-from-Kh}. Hence, one can use~\eqref{eq:s-equation} to deduce
the Rasmussen invariant of these knots from their homology.

There are only two knots with 16 crossings and $\GD=1$, for which the
assumption of~\ref{prop:s-from-Kh} fails. They are $16^n_{864894}$ and
$16^n_{925408}$. As it turns out, both of them satisfy the condition
of Lemma~\ref{lem:rank-difference} and, hence, \eqref{eq:s-equation} is still
applicable. Their Rasmussen invariants are both $0$.

\begin{table}[ht]
\tbl{All knots with at most 15 crossings and two knots with 16 crossings
that have Alexander polynomial 1 and non-zero Rasmussen invariant. The two
16-crossing knots or their mirror images can potentially be strongly
quasipositive.}{\normalsize%
\begin{tabular}{|l|c|c|c|c||l|c|c|c|c|}
\hline
knot $K$& $\!s(K)\!$& $\!hw(K)\!$& $\!e(K)\!$& $\!E(K)\!$&
knot $K$& $\!s(K)\!$& $\!hw(K)\!$& $\!e(K)\!$& $\!E(K)\!$\\
\hline\hline
$13^n_{1496}$  &  $2$& $3$&  $0$& $8$& $14^n_{7708}$  &  $2$& $3$&  $0$& $8$\\
\hline
$15^n_{28998}$ &  $2$& $3$& $-4$& $8$& $15^n_{87941}$ &  $2$& $3$& $-2$& $8$ \\
$15^n_{40132}$ &  $2$& $3$&  $0$& $8$& $15^n_{89822}$ & $-2$& $3$& $-8$& $0$ \\
$15^n_{52282}$ &  $2$& $3$&  $0$& $6$&
$\mathbf{15^n_{113775}}$\!&  $\mathbf2$& $3$&  $\mathbf2$& $12$\\
$15^n_{54221}$ & $-2$& $3$& $-8$& $2$& $15^n_{132396}$&  $2$& $3$&  $0$& $6$ \\
$15^n_{58433}$ &  $2$& $3$& $-2$& $8$& $15^n_{139256}$&  $2$& $3$&  $0$& $10$\\
$15^n_{58501}$ &  $2$& $3$& $-2$& $8$& $15^n_{145981}$& $-2$& $3$& $-8$& $0$ \\
$15^n_{65084}$ &  $2$& $3$&  $0$& $8$& $15^n_{165398}$&  $2$& $3$&  $0$& $6$ \\
$15^n_{65980}$ & $-2$& $3$& $-8$& $2$&               &     &    &     &     \\
\hline
$\mathbf{16^n_{412372}}$\!& $\mathbf{-2}$& $3$&$-12$&$\mathbf{-2}$&
$\mathbf{16^n_{955859}}$\!&  $\mathbf2$& $3$&  $\mathbf2$&$12$ \\
\hline
\end{tabular}
}
\label{tbl:non-slice}
\end{table}

In total, 82 knots with up to 16 crossings have Alexander polynomial 1
and non-zero Rasmussen invariant. Theorem~\ref{thm:Rasmussen-main} implies the
following.
\begin{proposition}\label{prop:not-slice}
All knots with up to 15 crossings from Table~\ref{tbl:non-slice} as well as 65
knots with 16 crossings are topologically locally-flatly slice, but not
smoothly slice.
\end{proposition}

These 65 knots with 16 crossings have table numbers 2601, 4787, 10734, 15919,
35456, 38567, 54888, 55405, 63905, 64312, 85435, 88272, 95001, 100099, 146445,
196836, 201101, 205822, 211749, 213930, 225414, 231486, 233317, 247683,
247710, 249903, 253331, 271353, 281590, 287865, 322069, 345376, 355871,
359271, 367431, 380325, 383790, 412372, 418128, 432810, 446116, 464148,
470729, 487352, 499458, 528093, 538818, 542632, 548142, 591990, 596477,
637428, 644951, 696243, 707728, 738706, 740175, 762813, 762859, 768960,
809057, 817682, 842714, 884475, and 955859.

At least three of these examples were already known to be non-slice. The knots
$13^n_{1496}$ and $16^n_{955859}$ were considered by Alexander
Stoimenow in~\cite[Example~1\footnote{Stoimenow arranges all the knots with a
given number of crossings into a single list. Hence, the knot $13^n_{1496}$ has
number $1496+4878=6374$ in~\cite{Stoimenow-example}, where 4878 is the number
of alternating knots with 13 crossings.}]{Stoimenow-example}
and~\cite[Example~4.1]{Stoimenow-sliceness}. The knot $15^n_{113775}$ is the
$(-3,5,7)$-pretzel knot. It was shown to be non-slice by
Rudolph~\cite{Rudolph-quasipositivity}. Another of his examples, the untwisted
positive double of the trefoil, has 19 crossings. Hence, it is not in the
table.

In general, all knots that are either strongly quasipositive or the mirror
image of a strongly quasipositive knot are not slice~\cite{Rudolph-sliceness}.
Let us show that only two more knots from Proposition~\ref{prop:not-slice}
can (potentially) fall into this category.

Let $P_K(v,z)$ be a HOMFLYPT polynomial of $K$, that is, the polynomial in
variables $v$ and $z$ determined by the skein relations
\begin{equation}\label{eq:HOMFLYPT-skein}
v^{-1}P_{\includegraphics[scale=0.45]{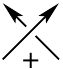}}
-vP_{\includegraphics[scale=0.45]{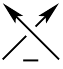}}
=zP_{\includegraphics[scale=0.45]{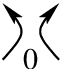}};
\qquad
P_{\includegraphics[scale=0.45]{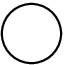}}=1.
\end{equation}
Let $e(K)$ and $E(K)$ be the lowest and highest exponents of $v$ in $P_K$,
respectively.

\begin{lemma}[{\cite[Theorem 1]{Morton-HOMFLYPT}}]\label{lem:Morton-degree}
Let $K$ be a knot represented by a diagram $D$. Then
\begin{equation}
w(D)-O(D)+1\le e(K)\le E(K)\le w(D)+O(D)-1.
\end{equation}
\end{lemma}

\begin{lemma}[cf.~\cite{Stoimenow-example}, Lemma 1]
1.~If $K$ is (strongly) quasipositive, then
\begin{equation}
0\le s(K)=2g_s(K)\le e(K).
\end{equation}
\noindent
2.~If the mirror image $\overline{K}$ of $K$ is (strongly) quasipositive, then
\begin{equation}
E(K)\le s(K)\le 0.
\end{equation}
\end{lemma}
\begin{proof}
Part 1 follows from~\ref{prop:s-for-quasipositive} and~\ref{lem:Morton-degree}
immediately. For part 2 notice, that $w(\overline D)=-w(D)$, $O(\overline
D)=O(D)$, where $\overline D$ is the mirror image of $D$, and that
$s(\overline K)=-s(K)$ (see~\cite[Proposition~3.11]{Rasmussen-s}).
\end{proof}

\begin{corollary}
Among all knots from Proposition~\ref{prop:not-slice}, only knots
$15^n_{113775}$, $16^n_{955859}$, and $\overline{16}^n_{412372}$ can be
strongly quasipositive (see Table~\ref{tbl:non-slice}). $15^n_{113775}$ is the
$(-3,5,7)$-pretzel knot and is indeed strongly quasipositive (see
Figure~\ref{fig:3-5-7-pretzel}).
\end{corollary}

\begin{remark}
After the original version of this paper was published, Stoimenow
demonstrated~\cite{Stoimenow-quasipositive} that neither $16^n_{955859}$ nor
$\overline{16}^n_{412372}$ are in fact strongly quasipositive.
\end{remark}

\subsection*{Acknowledgments}
I am very grateful to Mikhail Khovanov for bringing to my attention the remark
by Justin Roberts and Peter Teichner that the Rasmussen invariant shows that
certain knots with trivial Alexander polynomial are not smoothly
slice~\cite{Roberts-Teichner}. I started to work on this paper in an attempt
to understand their remark. I also thank Sebastian Baader for teaching me
quasipositive braids and knots.

\def\bibbox{\leavevmode\hbox}
\raggedright

\end{document}